\documentclass{article}
\usepackage{amsmath,amsthm,amssymb}
\usepackage{graphicx}

\newcommand{\CC}{\mathbb{C}}

\newcommand{\NN}{\mathbb{N}}

\title{\bf Mixed Neumann-Fourier expressions for solving integral equations} 
\author{Raimundas Vidunas\\
\em Vilnius University}

\begin{document}

\date{}
\maketitle

\begin{abstract}
While teaching a course on integral equations, I noticed that a straightforward combination of Neumann series 
and Fourier series for the resolvent (or the solution) of an integral equation has good approximation qualities. 
This short article presents and investigates this combination of approximating series.
\end{abstract}

\section{The mixed series}

Fredhol'ms integral equation of the second kind
\begin{equation} \label{eq:fredholm}
y(x)=f(x)+\lambda\int_a^b K(x,t) y(t)dt
\end{equation}
with $K(x,t)\in L_2([a,b]^2)$ can be solved \cite[Ch.~3]{Kanwal} in the resolvent form
\begin{equation} \label{fr:sol}
y(x)=f(x)+\lambda\int_a^b R(x,t,\lambda) f(t)dt.
\end{equation}
The resolvent $R(x,t,\lambda)$ can be expressed as the Neumann series
\begin{equation}
R(x,t,\lambda)=\sum_{n=1}^{\infty} \lambda^{n-1} K_n(x,t) 
\end{equation}
in terms of the iterated kernels:
\begin{align}
K_1(x,t)= &\; K(x,t), \\  \label{eq:iterk}
K_n(x,t)= &\; \int_a^b K(x,\xi)K_{n-1}(\xi,t)d\xi.
\end{align}
The resolvent is a meromorphic function of $\lambda\in\CC$,
with only $\lambda=\infty$ as a possible limiting point. 
The poles are the eigenvalues of the corresponding homogeneous equation
\begin{equation} \label{eq:homog}
\psi(x)=\lambda\int_a^b K(x,t) \psi(t)dt.
\end{equation}
The Neumann series converges for $|\lambda|<1/||K||$,  or up to the nearest pole. 
The integral equation (\ref{eq:fredholm}) has a unique solution when $\lambda$ is not an eigenvalue \cite[Ch.~4.2]{Kanwal}. 

If the kernel $K(x,t)\in L_2([a,b]^2)$ is symmetric, 
that is, if $K(x,t)=\overline{K(t,x)}$, then we have the Fourier series expression 
\begin{equation} \label{eq:rfourier}
R(x,t,\lambda)=K(x,t)+\lambda \,
\sum_{k=1}^{\infty} \frac{\psi_k(x)\,\overline{\psi_k(t)}}{\lambda_k\,(\lambda_k-\lambda)}
\end{equation}
in terms of an orthonormal system $\psi_1(x),\psi_2(x),\ldots$ of eigenfunctions \cite[Ch.~7]{Kanwal}
with the corresponding eigenvalues $\lambda_1,\lambda_2,\ldots$
of the homogeneous Fredholm equation (\ref{eq:homog}).
The Fourier series converges absolutely and almost uniformly
when $\lambda$ is not an eigenvalue \cite[p.~115]{Tricomi}.
And we have the bilinear series
\begin{equation}
K(x,t)=
\sum_{k=1}^{\infty} \frac{\psi_k(x)\,\overline{\psi_k(t)}}{\lambda_k}
\end{equation}
that converges in the $L_2$-norm generally. 
Accordingly, we can modify (\ref{eq:rfourier}) to
\begin{equation} \label{eq:rfourie}
R(x,t,\lambda)= 
\sum_{k=1}^{\infty} \frac{\psi_k(x)\,\overline{\psi_k(t)}}{\lambda_k-\lambda},
\end{equation}
but the convergence is weaker. 

On the other hand, we can improve convergence of (\ref{eq:rfourier})
by subtracting known Fourier series whose terms behave asymptotically as $\sim 1/\lambda_k^2$.
Suitable Fourier series are for the iterated kernels:
\begin{equation} 
K_n(x,t)=\sum_{k=1}^{\infty} \frac{\psi_k(x)\,\overline{\psi_k(t)}}{\lambda_k^n}.
\end{equation} 
For $n\ge 2$, these series converges absolutely and almost uniformly. 
Employing the series with $n=2$, we obtain
\begin{equation} \label{eq:rfourier1}
R(x,t,\lambda)=K(x,t)+\lambda K_2(x,t)+\lambda^2 
\sum_{k=1}^{\infty} \frac{\psi_k(x)\,\overline{\psi_k(t)}}{\lambda_k^2\,(\lambda_k-\lambda)},
\end{equation}
with the Fourier terms behaving asymptotically as $\sim 1/\lambda_k^3$.
Next, employing $K_3(x,t)$ we obtain
\[
R(x,t,\lambda)=K(x,t)+\lambda K_2(x,t)+\lambda^2 K_3(x,t)+\lambda^3
\sum_{k=1}^{\infty} \frac{\psi_k(x)\,\overline{\psi_k(t)}}{\lambda_k^3\,(\lambda_k-\lambda)}.
\]
Generally, we have
\begin{equation} \label{eq:rfourier2}
R(x,t,\lambda)=\sum_{n=1}^m \lambda^{n-1}K_n(x,t)+\lambda^m 
\sum_{k=1}^{\infty} \frac{\psi_k(x)\,\overline{\psi_k(t)}}{\lambda_k^m\,(\lambda_k-\lambda)}.
\end{equation}
The obvious practical advantage is that convergence of the Fourier series improves 
with increased $m$. The resolvent can be effectively approximated by a few terms of the Neumann series 
plus several terms of Fourier series. Even if $\lambda$ is near some large eigenvalue $\lambda_n$, 
the terms with $|\lambda_k|>|\lambda_n|$ diminish quickly. 

From the probed literature it appears that the closest topic that may suggest formulas (\ref{eq:rfourier1}), (\ref{eq:rfourier2})
is ``acceleration of convergence of Fourier series" or ``overcoming the Gibbs phenomenon" by means of {\em polynomial subtraction}
\cite[\S 5.1]{Adcock}, \cite{Nersesian}. The Neumann series part in (\ref{eq:rfourier2}) is often polynomial,
as the example of the following section illustrates.
The term ``Fourier-Neumann series" has been used to refer to 
Fourier series  in terms of Bessel functions \cite{Ciauri}.

\section{An example}

Consider the integral equation 
\begin{equation} \label{eq:inteq}
y(x)=x+\lambda\int_0^1 \min(x,t) y(t)dt,
\end{equation}
as in \cite[p.~124--127]{Kleiza}. 
The eigenvalues and eigenfunctions are found by solving the homogeneous equation
\begin{equation}
\psi(x)=\lambda\int_0^1 \min(x,t) \psi(t)dt,
\end{equation}
or the equivalent differential equation with boundary conditions:
\begin{equation} \label{eq:boundcond}
\psi''(x)+\lambda\psi(x)=0, \qquad \psi(0)=0, \qquad \psi'(1)=0.
\end{equation}
The  eigenvalues and the corresponding normalized eigenfunctions are parametrized by $n=1,2,\ldots$:
\begin{equation}
\lambda_n=\left(n-{\textstyle\frac12}\right)^2\pi^2,  \qquad
\psi_n=\sqrt2\,\sin \!\left(n-{\textstyle\frac12}\right)\!\pi x.
\end{equation}
By (\ref{eq:rfourie}), the resolvent is
\begin{equation} 
R(x,t,\lambda)= 
2 \, \sum_{n=1}^{\infty} \frac{\sin\!\left(n-\frac12\right)\!\pi x\;\sin\!\left(n-\frac12\right)\!\pi t}
{\left(n-\frac12\right)^2\pi^2-\lambda}.
\end{equation}
By (\ref{fr:sol}), we have
\begin{equation} \label{eq:approx1}
y(x)=x+
\lambda\,\sum_{n=1}^{\infty} \frac{2(-1)^{n+1}\sin\!\left(n-\frac12\right)\!\pi x}
{\left(n-\frac12\right)^2\!\pi^2\left(\!\left(n-\frac12\right)^2\!\pi^2-\lambda\right)}.
\end{equation}
To compute the iterated kernels for (\ref{eq:rfourier2}), 
let $\widetilde{K}_n(x,t)$ denote the restriction of $K_n(x,t)$ to $t\le x$. 
The iterated kernels are symmetric just as $K(x,t)$, and
\begin{align}
\widetilde{K}_{n+1}(x,t)= 
\int_0^t \xi\,\widetilde{K}_{n}(t,\xi)d\xi+\int_t^x \xi\,\widetilde{K}_{n}(\xi,t)d\xi
+x \int_x^1 \widetilde{K}_{n}(\xi,t)d\xi
\end{align}
by (\ref{eq:iterk}). Starting from $\widetilde{K}_1(x,t)=t$, we compute
\begin{align} 
\widetilde{K}_2(x,t)= &\; \textstyle  xt-\frac12\,x^2t-\frac16\,t^3,\\
\widetilde{K}_3(x,t)= &\; \textstyle  \frac13\,xt-\frac16\,x^3t-\frac16\,xt^3
+\frac1{24}\,{x}^{4}t+\frac{1}{12}\,{x}^{2}{t}^{3}+\frac1{120}\,t^5.
\end{align}
We compute the initial integrals 
\begin{equation} 
J_n(x)=\int_0^1 K_n(x,t)f(t)dt = \int_0^x \widetilde{K}_n(x,t)f(t)dt+ \int_x^1 \widetilde{K}_n(t,x)f(t)dt
\end{equation} 
with $f(x)=x$:
\begin{align} 
J_1(x)= &\; \textstyle  \frac12\,x-\frac16\,x^3,\\
J_2(x)= &\; \textstyle  \frac5{24}\,x-\frac1{12}\,x^3+\frac1{120}\,x^5,\\
J_3(x)= &\; \textstyle  \frac{61}{720}\,x-\frac{5}{144}\,x^3+\frac1{240}\,x^5-\frac{1}{5040}\,x^7.
\end{align}
We get the solution in subsequently quicker approximations:
\begin{align}   \label{eq:approx2}
y(x)= &\; x+\lambda J_1(x)+\lambda^2\,\sum_{n=1}^{\infty} \frac{2(-1)^{n+1}\sin\!\left(n-\frac12\right)\!\pi x}
{\left(n-\frac12\right)^4\!\pi^4\left(\!\left(n-\frac12\right)^2\!\pi^2-\lambda\right)},\\
 \label{eq:approx3}
y(x)= &\; x+\lambda J_1(x)+\lambda^2J_2(x)+\lambda^3\,
\sum_{n=1}^{\infty} \frac{2(-1)^{n+1}\sin\!\left(n-\frac12\right)\!\pi x}
{\left(n-\frac12\right)^6\!\pi^6\left(\!\left(n-\frac12\right)^2\!\pi^2-\lambda\right)},\\
 \label{eq:approx4}
y(x)= &\; x+\lambda J_1(x)+\lambda^2J_2(x)+\lambda^3J_3(x) \nonumber \\
& \, \hspace{95pt} + \lambda^4\,
\sum_{n=1}^{\infty} \frac{2(-1)^{n+1}\sin\!\left(n-\frac12\right)\!\pi x}
{\left(n-\frac12\right)^8\!\pi^8\left(\!\left(n-\frac12\right)^2\!\pi^2-\lambda\right)}.
\end{align}
The solution exists only if and only if $\lambda$ is not an eigenvalue.
For $\lambda>0$ or $\lambda<0$, an elementary expression for the solution is,
respectively,
\begin{align}
y(x)=\frac{\sin \sqrt{\lambda}\,x}{\sqrt{\lambda}\,\cos\sqrt\lambda} \qquad\mbox{or}\qquad
y(x)=\frac{\sinh \sqrt{-\lambda}\,x}{\sqrt{-\lambda}\,\cosh\sqrt{-\lambda}}.
\end{align}
The poles of these solutions at the eigenvalues $\lambda$ are well captured by the respective Fourier terms.
With those Fourier terms, the approximation error can be seen as a continuous function of $\lambda$. 

\begin{figure}
\begin{picture}(300,174)
\put(-9,0){\includegraphics[height=172pt]{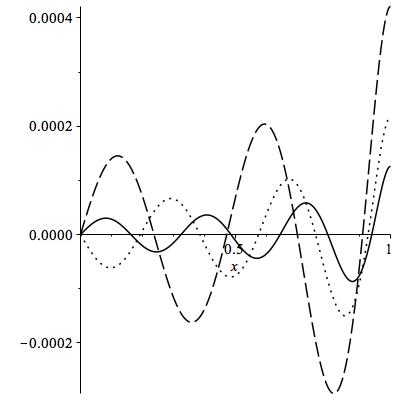}}
\put(175,0){\includegraphics[height=172pt]{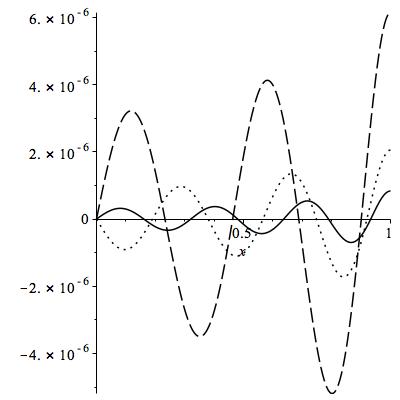}}
\end{picture}
\caption{Approximation error for the truncated series in (\ref{eq:approx1}) and (\ref{eq:approx2})
with 4 (dashed line), 5 (dotted line) or 6 (solid line) Fourier terms,
for $\lambda=4$, $x\in [0,1]$.}
  \label{fg:approx}
\end{figure}

With explicit solutions at hand, we may analyze approximation quality of differently 
truncated mixed Neumann-Fourier series (\ref{eq:approx1}), (\ref{eq:approx2})--(\ref{eq:approx4}).
Figure \ref{fg:approx} and Table \ref{tb:approx} represent a batch of relevant data. 
Figure \ref{fg:approx} demonstrates how adding just one Neumann series term diminishes
the approximation error by a factor $\approx 100$. The approximation error is largest at the endpoint $x=1$,
because the non-homogeneous term $f(x)=x$ in (\ref{eq:inteq}) does not satisfy the boundary condition (\ref{eq:boundcond})
at $x=1$, resulting in a weak version of the Gibbs phenomenon \cite[\S 4.2]{Adcock}. 
Table \ref{tb:approx} suggests that adding a Neumann series term is more worthwhile than adding a Fourier term,
unless $\lambda$ becomes large while the number of Fourier terms lags the number of Neumann terms.
The approximation errors for negative $\lambda$ are of the same order as respectively those for $|\lambda|$.
 
\begin{table}
\begin{tabular}{|c|c|llllll|}
\hline
 & $n$ & $\lambda=2$ & $\lambda=5$ & $\lambda=10$ & $\lambda=20$ & $\lambda=50$ & $\lambda=100$ \\
\hline
(\ref{eq:approx1}) & 3 & $4.9\;10^{-4}$ & $1.2\;10^{-3}$ & $2.6\;10^{-3}$ & $5.5\;10^{-3}$ & $1.8\;10^{-2}$ & $9.6\;10^{-2}$ \\
& 4 & $2.1\;10^{-4}$ & $5.3\;10^{-4}$ & $1.1\;10^{-3}$ & $2.2\;10^{-3}$ & $6.4\;10^{-3}$ & $1.7\;10^{-2}$ \\
& 5 & $1.1\;10^{-4}$ & $2.7\;10^{-4}$ & $5.5\;10^{-4}$ & $1.1\;10^{-3}$ & $3.1\;10^{-3}$ & $7.1\;10^{-3}$ \\
& 6 & $6.3\;10^{-5}$ & $1.6\;10^{-4}$ & $3.2\;10^{-4}$ & $6.5\;10^{-4}$ & $1.7\;10^{-3}$ & $3.8\;10^{-3}$ \\
\hline
(\ref{eq:approx2}) & 3 & $6.1\;10^{-6}$ & $3.9\;10^{-5}$ & $1.6\;10^{-4}$ & $7.1\;10^{-4}$ & $6.0\;10^{-3}$ & $7.2\;10^{-2}$ \\
& 4 & $1.5\;10^{-6}$ & $9.6\;10^{-6}$ & $3.9\;10^{-5}$ & $1.7\;10^{-4}$ & $1.2\;10^{-3}$ & $6.8\;10^{-3}$ \\
& 5 & $5.1\;10^{-7}$ & $3.2\;10^{-6}$ & $1.3\;10^{-5}$ & $5.4\;10^{-5}$ & $3.7\;10^{-4}$ & $1.8\;10^{-3}$ \\
& 6 & $2.1\;10^{-7}$ & $1.3\;10^{-6}$ & $5.3\;10^{-6}$ & $2.2\;10^{-5}$ & $1.4\;10^{-4}$ & $6.5\;10^{-4}$ \\
\hline
(\ref{eq:approx3}) & 3 & $8.9\;10^{-8}$ & $1.4\;10^{-6}$ & $1.2\;10^{-5}$ & $1.0\;10^{-4}$ & $2.3\;10^{-3}$ & $5.7\;10^{-2}$ \\
& 4 & $1.3\;10^{-8}$ & $2.1\;10^{-7}$ & $1.7\;10^{-6}$ & $1.4\;10^{-5}$ & $2.6\;10^{-4}$ & $3.0\;10^{-3}$ \\
& 5 & $2.8\;10^{-9}$ & $4.5\;10^{-8}$ & $3.6\;10^{-7}$ & $3.0\;10^{-6}$ & $5.2\;10^{-5}$ & $5.1\;10^{-4}$ \\
& 6 & $8.1\;10^{-10}$ & $1.3\;10^{-8}$ & $1.0\;10^{-7}$ & $8.4\;10^{-7}$ & $1.4\;10^{-5}$ & $1.3\;10^{-4}$ \\
\hline
(\ref{eq:approx4}) & 3 & $1.4\;10^{-9}$ & $5.5\;10^{-8}$ & $9.2\;10^{-7}$ & $1.6\;10^{-5}$ & $8.9\;10^{-4}$ & $4.6\;10^{-2}$ \\
& 4 & $1.2\;10^{-10}$ & $4.7\;10^{-9}$ & $7.7\;10^{-8}$ & $1.3\;10^{-6}$ & $6.0\;10^{-5}$ & $1.4\;10^{-3}$ \\
& 5 & $1.7\;10^{-11}$ & $6.7\;10^{-10}$ & $1.1\;10^{-8}$ & $1.8\;10^{-7}$ & $7.8\;10^{-6}$ & $1.5\;10^{-4}$ \\
& 6 & $3.4\;10^{-12}$ & $1.3\;10^{-10}$ & $2.2\;10^{-9}$ & $3.6\;10^{-8}$ & $1.5\;10^{-6}$ & $2.7\;10^{-5}$ \\
\hline
\end{tabular}
\caption{Approximation errors at $x=1$  for the truncated series (\ref{eq:approx1}), (\ref{eq:approx2})--(\ref{eq:approx4})
with $n\in\{3,4,5,6\}$ Fourier terms.}
 \label{tb:approx}
\end{table}

\section{More formulas}

Let us consider the interplay of Fourier series with Fredholm's representation 
\begin{equation}
R(x,t,\lambda)=\frac{D(x,t,\lambda)}{\Delta(\lambda)},
\end{equation}
where 
\begin{equation} \label{eq:delta}
\Delta(\lambda)=\prod_{k=1}^{\infty} \left(1-\frac{\lambda}{\lambda_k}\right)
\end{equation}
is an analytic function with the zeroes exactly at the eigenvalues \cite[Ch.~4]{Kanwal}.
The analytic series in $\lambda$ of $\Delta(x)$ and $D(x,t,\lambda)$  converge on the whole $\CC$
when the kernel $K(x,t)$ is $L_2$-integrable.
The results are unremarkable for approximation purposes, but there is some representational appeal.

The power series coefficients of 
\begin{align}
\Delta(\lambda) =\sum_{n=0}^{\infty} c_n\lambda^n, \qquad
D(x,t,\lambda)= \sum_{n=0}^{\infty} C_n(x,t)\lambda^n,
\end{align} 
can be computed \cite[(4.2.9)--(4.2.13a)]{Kanwal} from the recurrence relations
\begin{align}
c_n = &\; -\frac{1}{n}\int_a^b C_{n-1}(x,x)dx,\\
C_n(x,t)= &\; c_nK(x,t)+\int_a^b K(x,\xi)C_{n-1}(\xi,t)d\xi,
\end{align} 
and the initial conditions $c_0=1$, $C_0(x,t)=K(x,t)$. 
Keeping in mind orthonormality of the $\psi_k(\xi)$'s, manipulation of the Fourier series gives
\begin{align} \label{eq:sc1}
c_1 = &\,  -\sum_{k=1}^{\infty} \frac{1}{\lambda_k}, \\
C_1(x,t)=  &\, -\sum_{k=1}^{\infty} \sum_{j\neq k}  \frac{\psi_j(x)\overline{\psi_j(t)}}{\lambda_k\lambda_j} \\
=&\, -\mathop{\sum\sum}_{j<k} 
\frac{\psi_j(x)\overline{\psi_j(t)}+\psi_k(t)\overline{\psi_k(t)}}{\lambda_j\lambda_k}, \\ 
c_2 = & \;  \mathop{\sum\sum}_{j<k}  \frac{1}{\lambda_j\lambda_k}, \\
C_2(x,t)=   &\;  \mathop{\sum\sum}_{j<k} \! \sum_{i\not\in\{j,k\}} 
\frac{\psi_i(x)\overline{\psi_i(t)}}{\lambda_j\lambda_k\lambda_i} \\
= &\, \mathop{\sum\sum\sum}_{i<j<k} 
\frac{\psi_i(x)\overline{\psi_i(t)}+\psi_j(t)\overline{\psi_j(t)}+\psi_k(t)\overline{\psi_k(t)}}{\lambda_i\lambda_j\lambda_k}, 
\\ 
c_3 = & \, - \mathop{\sum\sum\sum}_{i<j<k}  \frac{1}{\lambda_i\lambda_j\lambda_k},
\end{align} 
etc. The pattern
\begin{align}
c_n = &\,  (-1)^n \! \sum_{S\subset \NN\;\&\,|S|=n} \Big(\prod_{j\in S} \, \lambda_j\Big)^{-1}, \\
C_n(x,t)=  &\, (-1)^n \! \sum_{S\subset \NN\;\&\,|S|=n} \Big(\prod_{j\in S} \, \lambda_j\Big)^{-1}
\, \sum_{j\in S} \psi_j(x)\overline{\psi_j(t)}.
\end{align} 
is apparent from (\ref{eq:delta}) and eventually from $D(x,t,\lambda)=\Delta(\lambda)R(x,t,\lambda)$ as well.
Towards an analogy to (\ref{eq:rfourier2}), we have
\begin{equation} \label{eq:fredhr}
D(x,t,\lambda)=\sum_{n=0}^{m-1} \lambda^{n}C_n(x,t)+ \!\!
\sum_{S\subset\NN\;\&\,|S|\ge m} \frac{(-\lambda)^{|S|}}{\prod_{j\in S} \lambda_j} \, \sum_{j\in S} \psi_j(x)\overline{\psi_j(t)}.
\end{equation}
But convergence is expected to be poor generally.
If the series (\ref{eq:sc1}) for $c_1$ does not converge, the series for the other $c_n$ 
should not converge (absolutely at least) either.
The series for $C_n(x,t)$ should then converge only in the $L_2$-norm. 

As an intermediate alternative to (\ref{eq:rfourier2}) and (\ref{eq:fredhr}) one may consider
\begin{align} \label{eq:pfourier}
\left(1-\alpha\lambda\right) R(x,t,\lambda)= & \,
K(x,t)+\sum_{n=1}^{m-1} \lambda^{n} \big(K_{n+1}(x,t)-\alpha K_{n}(x,t) \big) \nonumber \\
 & +\lambda^m 
\sum_{k=1}^{\infty} \frac{(1-\alpha\lambda_k)\,\psi_k(x)\,\overline{\psi_k(t)}}{\lambda_k^m\,(\lambda_k-\lambda)}.
\end{align}
This modification could be useful for $\lambda$ near some large eigenvalue $\lambda_k$,
taking \mbox{$\alpha=1/\lambda_k$} and annihilating the $k$th term in the Fourier series. 
Similar expressions can be obtained for $R(x,t,\lambda)$ multiplied by several such factors $(1-\alpha_j\lambda)$
annihilating several peaking Fourier terms. But the Fourier terms in (\ref{eq:pfourier}) decrease as $O(1/\lambda_j^m)$
rather than $O(1/\lambda_j^{m+1})$ in (\ref{eq:rfourier2}). Eliminating a peaking Fourier term comes
at the cost of one added Neumann series term for the purpose of better convergence of the Fourier series.

\end{document}